\documentclass[12pt]{amsart}

\usepackage{latexsym}
\usepackage{amsmath}
\usepackage{amssymb}
\usepackage{amscd}
\usepackage{psfig}
\usepackage{multicol}

\newtheorem{theorem}{\bf Theorem}[section]

\newtheorem{proposition}[theorem]{\bf Proposition}
\newtheorem{lemma}[theorem]{\bf Lemma}

\def\R{\mathbb{R}}
\def\C{\mathbb{C}}

\def\l{\lambda}
\title[Conjugacy classes in $SU(2)$]{Products of conjugacy classes in $SU(2)$}

\date{\today}
\begin{document}
\begin{abstract}
We obtain a complete description of the conjugacy classes
$C_1,\ldots,C_n$ in $SU(2)$ with the property that $C_1\ldots
C_n=SU(2)$. The basic instrument is a  characterization of the
conjugacy classes $C_1,\ldots,C_{n+1}$ in $SU(2)$ with $C_1\ldots
C_{n+1}\ni I$, which generalizes a result of \cite{Je-We}.
\end{abstract}

\author[L. C. Jeffrey]{Lisa  C. Jeffrey}
\address[L.\ C.\ Jeffrey]{Department of Mathematics\\ University of Toronto
 \\Toronto, Ontario M5S 3G3, Canada}
 \email{jeffrey@math.toronto.edu}

\author[A.-L. Mare]{Augustin-Liviu  Mare}
\address[A.-L.\ Mare]{Department of Mathematics\\ University of Toronto
 \\Toronto, Ontario M5S 3G3, Canada}
 \email{amare@math.toronto.edu}

\maketitle

\section{Introduction}

The following problem was posed by D. Burago:

{\bf Problem:} Let $G$ be a group. For which
conjugacy classes $C_1, \dots, C_n$ of $G$ is it true that
the multiplication map
$$ C_1 \times \dots \times C_n \to G$$
is surjective?

We give a solution to this problem in the case
$G = SU(2)$. In this case the conjugacy classes are
parametrized by their eigenvalues
$${\rm diag}(e^{i \lambda}, e^{-i \lambda} )$$
so they are determined by one number
$\lambda \in [0, \pi]$.

Burago's interest was primarily in discrete groups.
The purpose of this note is to point out that the problem he posed 
is also of interest for Lie groups such as $SU(2)$, and to exhibit
a solution in that case.

For more general Lie groups $G= SU(n)$ the problem could be 
studied by adapting results on the quantum cohomology 
of Grassmannians: see [Ag-Wo]. The problem is 
related to recent results described in the 
article [KLM]. 

\section{Eigenvalues of a multiple product}

For any $\l \in[0,\pi]$ we denote by $C(\l)$ the conjugacy class
of the matrix $${\rm diag}(e^{i\lambda},e^{-i\lambda})$$ in
$SU(2)$. Note that any conjugacy class in $SU(2)$ is of the form
$C(\lambda)$ for a unique $\lambda\in [0,\pi]$. The following
result was proved in [Je-We] (Proposition 3.1):

\begin{proposition}\label{prop:Je-We} For $\l_1,\l_2,\l_3\in [0,\pi]$ we have
$$C(\lambda_1)C(\lambda_2)C(\lambda_3)\ni I$$
iff \begin{equation} |\l_1 - \l_2| \le \l_3 \le {\rm
min}\{\l_1+\l_2, 2\pi -(\l_1+\l_2)\}
\end{equation}
\end{proposition}

Note that (1) is equivalent to
\begin{eqnarray*} \l_1 +\l_2+\l_3\le 2\pi \\
-\l_1 -\l_2 +\l_3 \le 0\\
 -\l_1 +\l_2 -\l_3 \le 0\\
 \l_1 -\l_2 -\l_3\le 0.\\
\end{eqnarray*}

 The goal of this section is to prove the more general result:
\begin{theorem}\label{thm:image} For $n\geq 2$ an integer and
$\l_1,\ldots,\l_{n+1}\in[0,\pi]$ we have
\begin{equation}\label{n+1}C(\l_1) \ldots C(\l_{n+1})\ni I\end{equation}
iff the following system of inequalities holds:

a) If $n+1=2k$ is an even number:
\begin{equation}\label{eqn:image-even} S_{n+1}^1(\{\l_i\}) \leq (n-1)\pi,\quad
S_{n+1}^3(\{\l_i\})\leq (n-3)\pi, \quad \ldots \quad S_{n+1}^{2k-1}(\{\l_i\})\le 0 \end{equation}
where $S_{n+1}^j(\{\l_i\})$ is any sum of the type $\sum_{i=1}^{n+1}\pm\l_i$
which contains exactly $j$ minus signs.

b) If $n+1=2k+1$ is an odd number:
\begin{equation}\label{eqn:image-odd} S_{n+1}^0(\{\l_i\})\le n\pi,\quad
S_{n+1}^2(\{\l_i\})\le (n-2)\pi,\quad \ldots\quad  S_{n+1}^{2k}(\{\l_i\})\le 0
\end{equation}
where $S_{n+1}^j(\{\l_i\})$ has the same meaning as before.
\end{theorem}

{\bf Remarks.}
1. A more concise way to express both (\ref{eqn:image-even}) and (\ref{eqn:image-odd})
 is
$$S_{n+1}^{n-2j}(\{\l_i\})\leq 2j\pi$$
for any $0\leq j \leq n/2$ and any sum of the type
$S_{n+1}^{n-2j}$.

2. An elementary computation involving the binomial formula
shows that the number of
inequalities in both (\ref{eqn:image-even}) and (\ref{eqn:image-odd}) is
$${{n+1}\choose 0}+{{n+1}\choose 2}+ \ldots ={{n+1}\choose 1}+{{n+1}\choose 3}+ \ldots = 2^n.$$

We will use induction on $n$ to prove this theorem. In order
to make the induction step we will need the following result:

\begin{lemma}\label{lemma} The condition (\ref{n+1}) holds iff there exists
$\l \in[0,\pi]$ such that
\begin{equation}\label{n}C(\l_1) \ldots C(\l_{n-1})C(\l)\ni I \end{equation} and
\begin{equation}\label{three}C(\lambda)C(\lambda_{n})C(\lambda_{n+1})\ni I.\end{equation}
\end{lemma}

\begin{proof}
The fundamental group of
the $(n+1)$-punctured sphere
$\Sigma_{n+1}$ in two dimensions
is the free group on $n$ generators, or the group
$$\Pi_n = \langle x_1, \dots, x_{n+1} ~|~x_1 \dots x_{n+1} = 1\rangle $$
with $n+1$ generators and one relation. We can form the
$(n+1)$-punctured sphere by gluing an $n$-punctured sphere and
a 3-punctured sphere along one of the boundary  components of
each. Call $S$ the common boundary resulting from this
construction and consider the fundamental groups of the two
components as follows:
$$\Pi_{n-1}=\langle x_1,\ldots,x_{n-1},x ~|~ x_1\ldots x_{n-1}x=1\rangle$$
and
$$\Pi_2=\langle x',x_n,x_{n+1} ~|~ x'x_nx_{n+1}=1\rangle,$$
where $x$ and $x'$ represent the loop  $S$ in each of the two components.
From the theorem of Seifert-van Kampen, we have that
\begin{equation}\label{pi}\Pi_{n}=(\Pi_{n-1}\times\Pi_2)/\langle xx'=1\rangle
\end{equation}

Now we consider representations of these groups into $G=SU(2)$.
The condition  (\ref{n+1}) is equivalent to the existence of a representation
 $\rho$ of $\Pi_{n+1}$ such that
$$\rho(x_i)\in C(\l_i)$$
for any $1\leq i \leq n+1$. From (\ref{pi}), this is equivalent to the existence
of a representation $\rho_{n-1}$ of $\Pi_{n-1}$ which coincides with
$\rho$ on $x_1,\ldots,x_{n-1}$,  and a representation $\rho_2$ of  $\Pi_2$ which
coincides with $\rho$ on  $x_{n}$ and $x_{n+1}$,
  and such that $\rho_{n-1}$ and $\rho_2$ satisfy
$$\rho_{n-1}(x)\rho_2(x')=I.$$
The latter equality implies that the conjugacy classes of $\rho_{n-1}(x)$
and $\rho_2(x')$ are equal, call them $C(\l)$ (note that in $SU(2)$ every element
is conjugate to its inverse). The conditions (\ref{n}) and (\ref{three})
correspond  respectively to the
representations $\rho_{n-1}$ and $\rho_2$.

  \end{proof}

{\it Proof of Theorem \ref{thm:image}} Just the induction step has
to be performed. We want to prove that
$$C(\l_1) \ldots C(\l_{n+1})\ni I$$
iff equation (\ref{eqn:image-even}) or (\ref{eqn:image-odd})
holds. Suppose that $n=2k$ is an even number. Condition (\ref{n})
 of Lemma \ref{lemma} is equivalent to
 \begin{equation}\label{all}  S_{n}^1(\l_1,\ldots ,\l_{n-1},\l) \leq (n-2)\pi,
S_{n}^3(\l_1,\ldots ,\l_{n-1},\l)\leq (n-4)\pi,  \ldots ,S_{n}^{2k-1}(\l_1,\ldots ,\l_{n-1},\l)\le 0
\end{equation}
where we have used the induction hypothesis, and  condition
(\ref{three})  is equivalent to
\begin{equation}\label{threen} |\l_n-\l_{n+1}|\leq \l \leq {\rm min}
\{\l_n+\l_{n+1}, 2\pi - (\l_n+\l_{n+1})\}\end{equation} where we
have used Proposition \ref{prop:Je-We}. By Lemma \ref{lemma},
condition (\ref{n+1}) is equivalent to the system of inequalities
obtained by considering each of the $2^{n-1}$ inequalities from
(\ref{all}) and deriving from it two inequalities, as follows:
\begin{itemize}
\item[(i)] if $\l$ occurs with a {\it plus } sign in that sum, replace it
 by  $\l_n-\l_{n+1}$ and $-\l_n+\l_{n+1}$
\item[(ii)] if $\l$ occurs with a {\it minus } sign in that sum, replace it
 by  $\l_n+\l_{n+1}$ and $-\l_n-\l_{n+1}$, but in the latter
situation add $2\pi$ to the right hand side of the original
inequality.
\end{itemize}
One sees that in the case (i) we replace an inequality of the  type
\begin{equation}\label{j}S_n^j\le (n-j  -1)\pi\end{equation} by two different
inequalities, both of the type
\begin{equation}\label{j+1}S_{n+1}^{j+1}\le (n-j-1)\pi.\end{equation}
In the case (ii) one again replaces an inequality of the type
(\ref{j}) by an inequality of the type (\ref{j+1}) and an
inequality of the type
$$S_{n+1}^{j-1}\leq (n-j+1)\pi.$$
One obtains $2^n$ distinct inequalities of type
(\ref{eqn:image-odd}), which means that (\ref{n+1}) is really
equivalent to (\ref{eqn:image-odd}).

A similar argument can be used when $n=2k-1$  is an odd number.
{\hfill{{$\square$}}\\}

{\bf Remark.} The result stated in Theorem 1.2 can also be
obtained from [Ag-Wo, Theorem 3.1] by using the structure of the
quantum cohomology ring of $\C P^1$. More precisely, let us
consider the two Schubert classes in $H^*(\C P^1)$:$$[\sigma_1]\in
H^2(\C P^1) \ {\rm and } \ [\sigma_2]=1\in H^0(\C P^1).$$ The
quantum cohomology ring of $\C P^1$ is
$$QH^*(\C P^1)=(H^*(\C P^1)\otimes \R[q],\star),$$
where $q$ is a formal variable of degree 4 and $\star$ is an
$\R[q]$-linear, commutative and associative   product which
satisfies
\begin{equation}\label{sigma}[\sigma_1]\star[\sigma_1]=q.\end{equation}
Each of the $2^n$ inequalities indicated in Theorem 2.2 can be
obtained by choosing  $i_1,\ldots, i_n\in\{1,2\}$ and evaluating
the product
$$[\sigma_{i_1}]\star \ldots \star [\sigma_{i_n}]$$
in $QH^*(\C P^1)$. By the equation (\ref{sigma}), this product is
of the form $q^d\sigma_{k},$ where $d$ is a positive integer and
$k\in\{1,2\}$.  The inequality of the type (\ref{eqn:image-even})
or (\ref{eqn:image-odd}) which corresponds to $i_1,\ldots, i_n$ is
$$\sum_{j=1}^n(-1)^{i_j-1}\lambda_j+(-1)^k\lambda_{n+1}
\leq 2d\pi.$$

\section{Surjectivity of a multiple product}

Our main result is
\begin{theorem} We have
\begin{equation}\label{surj}C(\l_1)\ldots C(\l_n)=SU(2)\end{equation}
iff for any integer  $j$ with  $0\le j < n/2$ and for any sum of
the type $S_n^j=S_n^j(\{\l_i\})$ (see Thm.1.2)  we have
\begin{equation}\label{equ}-(j-1)\pi \le S_n^j\le (n-j-1)\pi .\end{equation}
 \end{theorem}

\begin{proof} The idea of the proof is that (\ref{surj}) holds iff (\ref{n+1}) holds for
any $\l_{n+1}\in [0,\pi]$. In turn, (\ref{n+1}) is equivalent to
(\ref{eqn:image-even}) and (\ref{eqn:image-odd}). We
just have to take each inequality from (\ref{eqn:image-even})
(respectively (\ref{eqn:image-odd})) and make the following formal
replacements in its left hand side:
\begin{itemize}
\item[(i)]  $\lambda_{n+1}$ by $\pi$
\item[(ii)] $-\l_{n+1}$ by $0$.
\end{itemize}

Let us consider the case $n=2k-1$. We have to show that if we
perform (i) and (ii) for each inequality contained in (2), we
obtain exactly one of the following inequalities:
\begin{eqnarray} \label{zero}\pi\le S_n^0 \le (n-1)\pi\\
                   \label{unu}0\le S_n^1\le (n-2)\pi
                   \end{eqnarray}
\begin{eqnarray} \label{doi}-\pi \le S_n^2\le (n-3)\pi\\
                  \label{trei}-2\pi \le S_n^3\le (n-4)\pi
\end{eqnarray}
\begin{eqnarray*} ..............................\end{eqnarray*}
We claim that if we label the inequalities given by (2) as [1],
[3],$\ldots$, [$2k-3$], [$2k-1$] then each of [1] and [$2k-1$]
gives exactly one of (\ref{zero}) and (\ref{unu}), each of [3] and
[$2k-3$] gives exactly one of (\ref{doi}) and (\ref{trei}),
$\ldots$ and finally

- if $k=2p$ is even, then each of [$2p-1$] and [$2p+1$] gives
exactly one of
\begin{eqnarray*} -(k-1)\pi \le S_n^{k-2}\le (k-2)\pi\\
                 -(k-2)\pi \le S_n^{k-1}\le (k-1)\pi
                 \end{eqnarray*}

- if $k=2p+1$ is odd, then each of [$2p+1$] gives exactly one of
\begin{eqnarray*} -(k-2)\pi \le S_n^{k-1}\le (k-1)\pi
\end{eqnarray*}

Consider first [1] together with [$2k-1$]: the only $S_{n+1}^1$
which contains $-\l_{n+1}$ leads to
$$\l_1+\ldots + \l_n \le (n-1)\pi$$
whereas the only $S_{n+1}^{2k-1}$ which contains $\l_{n+1}$ leads
to
$$\l_1 +\ldots +\l_n\ge \pi.$$
The remaining inequalities of type $S_{n+1}^1\le (n-1)\pi$ lead to
all possible inequalities of the type
$$S_n^1 \le (n-2)\pi$$
and the remaining inequalities of the type $S_{n+1}^{2k-1}\le 0$
lead to all possible inequalities of the type
$$S_n^1 \ge 0.$$

The same idea applies\footnote{If we compare the total number of
inequalities we start with to the number of inequalities obtained
via (i) and (ii) we ``deduce" that $ {{n+1}\choose
{2j+1}}+{{n+1}\choose{n+1-(2j+1)}}=2({n\choose{2j+1}}+{n\choose{2j}}).$
The latter equation is obviously true, by  properties of 
Pascal's triangle.} to each pair [$2j+1$], [$2(k-j)-1$], $0\le j<
k/2$ (if $k=2p+1$ is an odd number, then for $j=p$ we have
$2j+1=2(k-j)-1$ and the corresponding pair reduces to just one
type of inequalities).

Similar ideas can be used in the case when $n=2k$ is an even
number.
\end{proof}

{\bf Remark.} The system of inequalities (\ref{equ}) admit
solutions for any $n\geq 2$. For $n=2$ the  {\it unique} solution
is
\begin{equation}\label{last}
\lambda_1=\lambda_2=\frac{\pi}{2}.\end{equation} For $n\geq 3$
there are several solutions, one of them consisting of
$\l_1,\l_2$ given by (\ref{last}) and $$\l_3= \ldots \l_n=0.$$

\bibliographystyle{abbrv}

\end{document}